\documentclass{amsart}
\usepackage{amssymb,latexsym}
\usepackage{amsmath}
\usepackage{graphics}

\newtheorem{theorem}{Theorem}[section]
\newtheorem{corollary}[theorem]{Corollary}
\newtheorem{lemma}[theorem]{Lemma}
\newtheorem{remark}[theorem]{Remark}
\newtheorem{prop}[theorem]{Proposition}
\newtheorem{definition}[theorem]{Definition}

\newtheorem{conjecture}[theorem]{Conjecture}
\newtheorem*{conj}{Graham's Conjecture}

\newcommand{\tPn}{\tilde P_n}

\newcommand{\tCn}{\tilde C_{n}}

\newcommand{\tD}{\tilde D}
\newcommand{\tH}{\tilde H}
\newcommand{\dist}{\text{dist}}
\begin{document}
\title{Cover pebbling cycles and certain graph products}
\author{Maggy Tomova}
\address{Department of Mathematics, University of California, Santa
Barbara, CA 93117} \email{maggy@math.ucsb.edu}

\author{Cindy Wyels}
\address{Department of Mathematics, California Lutheran University, 
Thousand Oaks, CA 91360} \email{wyels@clunet.edu}

\begin{abstract}
A pebbling step on a graph consists of removing two pebbles from one vertex
and placing one pebble on an adjacent vertex. A graph is said to be
{\it cover pebbled} if every vertex has a pebble on it after a series
of pebbling steps.
The \textit{cover pebbling number} of a graph is the minimum number
of pebbles such that the graph can be cover pebbled, no matter how
the pebbles are initially placed on the
vertices of the graph.
In this paper we determine the cover pebbling numbers of cycles, finite 
products of paths and cycles, and products of a path or a cycle with 
{\em good graphs}, amongst which are trees and complete graphs.
In the process we provide evidence in 
support of an affirmative answer to a question posed in a paper by Cundiff, Crull, et al.

\vspace{.1in}

\vspace{0.1in}

\end{abstract}

\maketitle

\section{Introduction}\label{S:intro}

The game of pebbling was first suggested by Lagarias and Saks as a tool for
solving a number-theoretical conjecture of Erd\"os. 
Chung successfully used this tool to prove the conjecture and 
established other 
results concerning pebbling numbers. In doing so she introduced pebbling 
to the literature \cite{C}.  

Begin with a graph $G$ and a certain number of pebbles placed on 
its vertices. A pebbling step consists of removing two pebbles from one vertex and 
placing
one pebble on an adjacent vertex. In (regular) pebbling, a target vertex is
selected, and the goal is to move a pebble to the target vertex. The 
minimum number of pebbles such that, regardless of their initial 
placement and regardless of the target vertex,  we can pebble that 
vertex is called {\em the pebbling number of $G$}.
In cover pebbling, the goal is to cover all the vertices with pebbles, 
i.e., to
move a pebble to every vertex of the graph simultaneously. 
The minimum number of pebbles required such that, regardless of their 
initial placement on $G$, there is a sequence of pebbling steps at the end of 
which every vertex has at least one pebble on it is called {\em the 
cover pebbling number of $G$}.
In the paper in 
which 
the concept of cover pebbling is introduced, the authors find
the cover pebbling numbers of several families of graphs, including trees 
and complete graphs 
\cite{H2}.
Hurlbert and Munyan have also announced a proof for the cover pebbling number of 
the
$n$-dimensional cube.

In this paper we ``translate" a distribution on a product of graphs 
to a distribution on one of the factors by introducing colors. This 
allows us to find upper bounds for the cover pebbling numbers of 
$G\square P_n$ (Corollary \ref{cor:path inequality}) and $G \square C_n$ 
(Corollary \ref{cor:cycle inequality}), 
where $G$ is any graph. As finding 
lower bounds given a particular graph is generally straightforward,  
in Corollary \ref{cor:cycle} we establish the cover pebbling number
of cycles. It is possible that 
upper bounds for the cover pebbling numbers of other products can be obtained using this technique.

Let $G = (V,E)$ be any graph.  A {\em distribution} of pebbles to the vertices 
of $G$ is any initial arrangement of pebbles on some subset $S$ of $V$. 
The set $S$ is called the {\em support} for the distribution;  vertices in $S$ are called support vertices.  
A {\em simple} distribution is one with a single support vertex. We use $\gamma (G)$ to denote the cover pebbling number of $G$. 

\begin{definition}
    A graph $G$ is {\em good} if 
    $$ \gamma (G) = \sum _{w \in V(G)} 2^{\dist (w,u)}$$
    for some vertex $u \in V(G)$.  Any vertex $u$ satisfying this 
    equation is a {\em key} vertex.
    \end{definition}

    \begin{remark}
\textnormal{A graph is good precisely when its cover pebbling number 
is equal 
to the number of pebbles needed to cover pebble the graph from a 
single (specific) vertex, i.e. from a key vertex. Thus when finding 
the cover pebbling number of a good graph, we only need to consider 
simple distributions.}
\end{remark}
In \cite{H2} we see 
that paths, trees and complete graphs are good, and the authors raise the question of 
whether every graph is good.  We believe this is the case:

\begin{conjecture}\label{conj:allgood}
    Every graph is good.
    \end{conjecture}
    
In support, we show that cycles are good. We also demonstrate that the 
product of any good graph with a cycle or a path is again good 
(Corollary \ref{cor:products of good}).

Chung's seminal pebbling result relies on products, and she lays out Graham's conjecture, 
perhaps the best known open question in pebbling.  
Say graphs $G$ and $H$ have vertex sets $V(G)=\{w_1, \dots , w_g\}$ and 
$V(H)=\{v_1, \dots , v_h\}$, respectively.  The product of $G$ and $H$, $G \square H$, is the graph 
with vertex set $V(G) \times V(H)$ (Cartesian product) and with 
edge set
\begin{align}
E(G\square H) = &\{\bigl( (w_1,v_1),(w_2,v_2) \bigr) \vert\, w_1=w_2 \text{ 
and } (v_1,v_2)\in V(H)\} \notag \\
&\cup \{\bigl( (w_1,v_1),(w_2,v_2) \bigr) \vert\, v_1=v_2 \text{ and } 
(w_1,w_2)\in V(G)\} \notag  .
\end{align}
Let $f(G)$ denote the pebbling number of the graph $G$.  
\begin{conj} $f(G \square H) \leq f(G)f(H)$.
    \end{conj}
    
There is much evidence in support of Graham's conjecture.  (See, for 
example, \cite{C}, \cite{E1}, and \cite{E2}.)  We believe that the analogous statement for 
cover pebbling 
involves equality:

\begin{conjecture}\label{conj:equality} $\gamma(G \square H)=\gamma(G)\gamma(H)$.
    \end{conjecture}

In Theorem 
\ref{thm:relationship} we show a relationship between Conjectures \ref{conj:allgood} 
and \ref{conj:equality} and in Lemma \ref{lem:equality} we demonstrate that Conjecture 
\ref{conj:equality} holds when $G$ is any good graph and $H$ is a path 
or a cycle. This allows us to easily compute the pebbling 
numbers of a large family of products, as
shown in Theorem \ref{thm:results}. In particular, we have proven the 
cover pebbling number of a finite product of cycles and paths which 
then yields the cover pebbling numbers of hypercubes ($P_2^n$), web graphs ($C_n \square P_m$), 
grids ($P_n \square P_m$), etc.

\section{Cover pebbling $G \square P_n$}

Let $G$ and $H$ be two graphs with vertices $w_1,..,w_g$ and 
$v_1,\ldots,v_h$ respectively. 
Recall that their product has vertex set 
$\{(w_i,v_j)|i=1,\ldots,g;\,  j=1,\ldots,h\}$.
We will associate to each distribution on $G \square H$ a certain 
distribution of colored pebbles on $H$. In some cases, namely when 
$H$ is a path or a cycle, results about this colored 
distribution can then be interpreted to obtain upper bounds for the 
pebbling number of the product.

We will call a distribution $t$-colored (or a $t$-distribution) if  each pebble in the 
distribution has been assigned one of $t$ possible colors. A color-respecting pebbling 
step for a colored distribution consists of taking two pebbles of the 
same color from some vertex and placing one of these pebbles on an adjacent vertex. 
When considering colored distributions we 
 allow only color-respecting steps. 
A distribution is $Q$-coverable if we can pebble the graph 
with $Q$ pebbles of any color on each vertex (performing only color-respecting steps). 
Thus the notions of cover pebbling and 
coverable distributions correspond to the case where $Q=t=1$.

To each distribution $D$ on $G \square H$ we associate a color 
distribution $\tD$ on $H$ in the following way:  use colors $c_1, c_2, \dots , c_g$ 
to assign color $c_i$ to each pebble that $D$ places on vertices $(w_i,v_j)$ (for any $j$).  
Collapse $G\square H$ to a single copy of $H$, which we call $\tH$ for clarity, by identifying 
$G \square \{v_i\}$ in $G \square H$ with vertex $V_i$ in $\tH$.  We place 
all pebbles from $G\square \{v_i\}$ on $V_i$.

\begin{lemma}\label{lem:equivalent}
 Let $G$ and $H$ be graphs and $D$  be a distribution on $G \square H$. 
 If the associated g-distribution $\tD$ on $\tH$ is $\gamma(G)$-coverable, then $D$ is 
 coverable.  
 \end{lemma}    

\begin{proof}
By hypothesis there is a sequence 
of color-respecting pebbling steps beginning with $\tD$ at the end of which there are 
$\gamma(G)$ pebbles 
on each vertex 
of $\tH$. Because the steps respect color we could have performed them in 
$G \square H$: taking two pebbles of color $c_i$ from $V_j$, discarding 
one and placing the other one on $V_k$ in 
$\tH$ corresponds to taking two pebbles from vertex $(w_i,v_j)$ 
and placing one of them on $(w_i,v_k)$. 
So there is a sequence of steps on $G \square H$, consisting 
only of moving pebbles from one copy of $G$ to another, at the end 
of which each copy of $G$ has $\gamma(G)$ pebbles. Now each copy of $G$ may be 
cover-pebbled using the $\gamma (G)$ pebbles on it, so $D$ is coverable.
\end{proof}

A priori it is possible that there exist coverable distributions on $G \square H$ 
that have associated colored distributions on $\tH$ that
are not $\gamma (G)$-coverable. However, in many cases it appears that
considering the color distribution on one of the 
factors is sufficient to find the pebbling number of the product. 

Within the usual concept of cover pebbling ($t$=1), given $M$ pebbles on a 
vertex $v$ we can always move $\lfloor \frac{M}{2}\rfloor$ pebbles to an 
adjacent vertex, possibly having to leave one pebble on $v$ in the case when $M$ is odd. 
The analogous statement holds for colored distributions.

\begin{lemma} \label{lem:movable}
  Suppose a vertex $v$ in the support of a t-colored distribution has $M>t$ 
  pebbles. Given any integer $E \leq M-t$, 
  at least $\lfloor E/2 \rfloor$ pebbles initially on $v$ can be 
  placed on an adjacent vertex using color-respecting steps.  

\end{lemma}    

\begin{proof}
Consider the set $T$ of all pebbles on the given vertex $v$. We will construct 
a subset $S$ of size at least $M-t$ consisting of pebbles all of which 
can be placed in same-color pairs. If a color has 
an odd number of representatives in $T$ remove one pebble of that color. As 
there are only $t$ colors, at most $t$ pebbles are removed. Let $S$ be 
the subset of all remaining pebbles, $|S| \geq M-t$. 
Now by removing pebbles in pairs of the same color we can obtain a 
smaller set, also containing even numbers of pebbles of each color, of 
size $E$ if $E$ is even or of size $E-1$ if $E$ is odd. Half of all pebbles in a given color 
can be moved to an adjacent vertex while discarding the other half. Thus we can move at least 
$\lfloor E/2\rfloor$ pebbles to an adjacent vertex. 
\end{proof}    

For the rest of this paper we will denote by $|V_s,\ldots,V_t|$ the number of pebbles on the path 
$V_s,\ldots,V_t$. The path on $m$ vertices will be denoted $P_m$.

The next proposition is slightly technical. 
The basic idea is that if we have a path on $m$ vertices and a 
distribution which places at least $Q$ pebbles on each of 
$V_2,\ldots,V_m$ and has $Q2^{m-1}$ additional pebbles, then we can use 
these 
additional pebbles to get at least $Q$ pebbles on $V_1$ and thus 
complete the $Q$-covering of the path. 

\begin{prop} \label{prop:cover v1}
Let $P_m$ be a path with at least 2 vertices, $Q$ and $K$ be integers with $Q>K$ and 
$Q \geq g$, and $\tD$ a g-distribution of $Q(m-1)+2^{m-1}Q$ pebbles such that 
$|V_i| \geq Q $ for all $i>1$ and $|V_1|=K$. Then there exists a sequence of color-respecting 
pebbling steps at the end of which $|V_i| \geq Q$ for all $1 \leq i \leq n$. 
\end{prop}

\begin{proof}
  We will use induction to prove this statement. If $m=2$, 
  $|V_2|=Q+2Q-K \geq Q+2Q-2K$. By Lemma \ref{lem:movable} with 
  $E=2Q-2K$ we can place $Q-K$ pebbles on $V_1$ leaving at least $Q$ 
  pebbles on $V_2$.  
 
  Now assume the result holds for all $i <m$. As $|V_1|=K$ and $|V_i|\geq Q$ for 
  $i=2,\dots ,m-1$, we know $|V_m| \leq Q + 2^{m-1}Q-K$.  Let $E=|V_m|-Q \leq 2^{m-1}Q-K$. 
  By Lemma \ref{lem:movable} we can leave $Q$ pebbles on $V_m$ while moving 
  $\lfloor E/2\rfloor$ pebbles from $V_m$ to $V_{m-1}$.  Prior to this move,
\begin{align*}
|V_1, \dots ,V_{m-1}| = &Q(m-1)+2^{m-1}Q-|V_m| \\
= &Q(m-1)+2^{m-1}Q-(E+Q)
\end{align*}
After moving $\lfloor E/2\rfloor$ pebbles to $V_{m-1}$ we have 
\begin{align*}
|V_1, \dots ,V_{m-1}| = &Q(m-1)+2^{m-1}Q-(E+Q)+ \lfloor E/2\rfloor \\
= & Q(m-2)+2^{m-1}Q-\lceil E/2\rceil \\
\geq & Q(m-2)+2^{m-1}Q-2^{m-2}Q \\
= &Q(m-2)+2^{m-2}Q.
\end{align*}
Thus by our induction hypothesis applied to the path 
$V_1,\ldots,V_{m-1}$ we can place $Q-K$ pebbles on $V_1$ for a total of 
$Q$ pebbles on $V_1$ while keeping at least $Q$ pebbles on each of $V_2,\ldots,V_m$. 
  \end{proof}

\begin{theorem}\label{thm:G times P}
     Let $g$ and $Q$ be positive integers with $g <Q$. If $\tD$ is any 
    $g$-distribution of $Q(2^n-1)$ pebbles on the vertices of $\tPn$, 
    then $\tD$ is Q-coverable.
 \end{theorem}   
 
 \begin{proof}

 Assume the theorem holds for all $m<n$. Let $\tD$ be a $g$-distribution on the vertices of 
$\tPn$. Label the vertices of $\tPn$ sequentially and so that
  $|V_1| \leq 
 |V_n|$ and let $K=|V_1|$.
 
 Case 1: $K \leq Q$
 
In this case $|V_2,\ldots,V_n|=Q(2^{n}-1)-K \geq 
 Q(2^{n}-1)-Q=Q(2^{n}-2) \geq Q(2^{n-1}-1)$. By the induction hypothesis we can 
 $Q$-cover the path $V_2,\ldots,V_n$ using at most $Q(2^{n-1}-1)$ pebbles. 
 Note that, as we only needed $Q(2^{n-1}-1)$ pebbles to $Q$-cover 
 $V_2,\ldots,V_{n}$, we now have $Q$ pebbles on each of   
 $V_2,\ldots,V_{n}$ and an additional 
 $Q(2^n-1)-K-Q(2^{n-1}-1)=Q2^{n-1}-K$ pebbles lying on the path $V_2,\ldots,V_{n}$. Thus 
 the path $V_1,\ldots,V_n$ now has a total of $Q(n-1)+Q(2^{n-1})$ 
 pebbles with at least $Q$ on each of $V_2,\ldots,V_n$ and $K$ 
 pebbles on $V_1$. 
 By Proposition \ref{prop:cover v1} we can move $Q-K$ pebbles to $V_1$ 
 keeping at least $Q$ pebbles at all other vertices, thus completing the 
 $Q$-covering of $\tPn$.
 
 Case 2: $K> Q$
 
Let $s$ be the largest integer 
such that for all $i \leq s$ the path $V_1,\ldots,V_i$ contains at 
least $Q(2^i-1)$ pebbles. Note that $s \geq 1$ since $K>Q$. By assumption $|V_n|\geq |V_1|>Q$
so $V_n$ is already $Q$-covered. If $s \geq n-1$ the 
 pebbles on $V_1,\ldots,V_{n-1}$ suffice to $Q$-cover $V_1,\ldots,V_{n-1}$ by 
the inductive hypothesis, so the distribution is $Q$-coverable and we are 
done. Thus we may assume $s \leq n-2$. $Q$-cover the path $V_1,\ldots,V_s$ 
using the pebbles lying on it (as is possible by the induction 
hypothesis). Note that $|V_1,\ldots,V_s| \leq Q(2^{s+1}-2)$ otherwise 
a larger integer $s$ could have been chosen. Now consider the path $V_{s+1},..,V_n$ 
which has $n-s$ 
vertices. It must have $Q(2^n-1)-|V_1,\ldots,V_s| > 
Q(2^n-1)-Q(2^{s+1}-1) =Q(2^n-1-2^{s+1}+1)=Q(2^n-2^{s+1}) \geq 
Q(2^n-2^{n-1})=Q(2^{n-1}) \geq Q(2^{n-s}-1) $ pebbles,  so by hypothesis we can 
$Q$-cover $V_{s+1},\ldots,V_n$. 

In either case $\tD$ is $Q$-coverable. 
\end{proof}

\begin{corollary} \label{cor:path inequality}
    For any graph $G$, $\gamma (P_n \square G) \leq\gamma(G)(2^n-1)$.
\end{corollary}
\begin{proof}
  Let $D$ be any distribution of $\gamma(G)(2^n-1)$ pebbles on $(P_n \square 
G)$. Letting $\gamma(G)=Q$ and $g$ be the number of vertices in $G$, 
by Theorem \ref{thm:G times P} we conclude that the associated 
$g$-distribution $\tD$ on $\tPn$ is $Q$-coverable. The result then follows from Lemma 
\ref{lem:equivalent}.
  
 \end{proof}   
 
 By letting $G$ consist of a single vertex we recover a result in \cite{H2}.
\begin{corollary}\label{cor:path}
   $\gamma(P_n)=2^n-1$ and $P_n$ is good.
 \end{corollary}
 \begin{proof}
    From Corollary \ref{cor:path inequality} we know $\gamma(P_n) \leq 2^n-1$. 
    To show $\gamma (P_n)\geq 2^n-1$, label the vertices of $P_n$ 
    sequentially and consider a distribution with 
    $v_1$ as the only support vertex. Covering $v_i$ from $v_1$ requires $2^{i-1}$, 
    pebbles so covering the whole path requires $\sum _{i=1}^{n}2^{i-1}=2^n-1$ pebbles.
  \end{proof}

\section {Cover Pebbling Number for Cycles}

In this section we obtain an upper bound for the cover pebbling number of the product of a cycle with any graph.  
A special case then gives the cover pebbling number of cycles.  Specifically, we show
$$\gamma(C_n \square G)\leq 
\begin{cases}
(2^{(n/2)+1}+2^{n/2}-3)\gamma(G), &\text{when $n$ is even;}\\
(2^{(n+1)/2}+2^{(n+1)/2}-3)\gamma(G), &\text{when $n$ is odd.}
\end{cases}
$$
In 
particular, taking $G$ to be a single vertex
$$\gamma(C_n)\leq 
\begin{cases}
2^{(n/2)+1}+2^{n/2}-3, &\text{when $n$ is even;}\\
2^{(n+1)/2}+2^{(n+1)/2}-3, &\text{when $n$ is odd.}
\end{cases}
$$

Fix some integer $n\geq 3$ and let $C_n$ be a cycle graph with vertices $V=\{v_1,...,v_n\}$, labeled 
sequentially. To simplify our discussion we let $r=n/2$ if $n$ is even and 
$r=(n+1)/2$ 
if $n$ is odd.  Let $P=2^r+2^{n-r+1}-3$; we will show that 
$\gamma(C_n)=P$. Let $G$ be any graph with $g$ vertices and let $\gamma(G)=Q$.

For the rest of this section we take $D$ be a distribution on $C_n \square 
 G$ and $\tD$ to be its associated $g$-distribution on $\tCn$. We will 
refer to a set $V_i,V_{i+1},\ldots,V_{r+i-1}$ of vertices of $\tCn$ as {\it primary}
when $V_i$ is a support 
vertex;  if $V_i$ is not necessarily a support vertex we will refer to the 
set as {\it secondary}.  Both primary and secondary sets are paths on $r$ 
vertices.  We 
will call a primary or secondary set {\it saturated} if it contains at 
least $Q(2^r-1)$ pebbles.

\begin{remark} \label {rmk:disjoint}
   \textnormal{Note that if, after 
color-respecting pebbling steps of the pebbles in $\tD$, there 
exists a partition of $\tCn$ into disjoint paths such that each of 
these paths has length $s_i$ and contains at least $Q(2^{s_i}-1)$ pebbles, 
then $\tD$ is $Q$-coverable by Theorem \ref{thm:G times P}}.
\end{remark}

\begin{lemma}\label{lem:one of each general}
    Suppose $D$ is a non-coverable distribution placing $PQ$ pebbles on 
$G \square C_n$. Let
    $\tD$ be its associated g-distribution on $\tCn$. 
    Then there is a sequential numbering of the vertices of $\tCn$ such that 
    $V_1,\ldots,V_r$ is saturated. With any such labeling there 
    exists  $i \leq r+1$ such that $V_i,\ldots,V_{r+i-1}$ is not 
    saturated.
      
\end{lemma} 
    
\begin{proof}
    By Lemma \ref{lem:equivalent} $\tD$ is not Q-coverable.
    
    {\bf Case 1: $n$ is even}
    
First consider the sets $V_1,\ldots,V_r$ and $V_{r+1},\ldots,V_n$. 
They cannot both be saturated otherwise $\tD$ would be coverable by Remark \ref{rmk:disjoint}. 
If both sets are unsaturated, then the total number of 
pebbles on the graph will be at most $Q((2^r-2)+(2^r-2)) 
<Q(2^{(n-r+1)}-1)+Q(2^{r}-2)=QP$,  
leading to a contradiction.  Thus one of the sets must be saturated and 
the other set must be unsaturated. After possibly relabeling the vertices 
$V'_1,\ldots,V'_n$ with $V'_1=V_{r+1}$ we have produced a labeling 
satisfying the conclusion of this lemma.

    {\bf Case 2: $n$ is odd}
    
First consider the sets $V_1,..,V_r$ and $V_{r+1},\ldots,V_n, V_1$. 
If both sets are unsaturated then we have at most 
$Q((2^r-2)+(2^r-2)) <QP$ pebbles in $\tD$, therefore one of them must be 
saturated. After possibly letting $V'_1=V_r$ we may assume 
$V'_1,\ldots,V'_r$ is saturated. If one of $V'_r,..,V'_n$ and 
$V'_{r+1},..,V'_{n}, V'_1$
is unsaturated we would be done, so suppose they are both saturated. 
By Remark \ref{rmk:disjoint} we can assume 
 $V'_{r+1},..,V'_{n}$ contains at most  $Q(2^{r-1}-2)$ pebbles as 
 $V_1,\ldots,V_r$ contains at least $Q(2^r-1)$ pebbles. That means $V'_1$ and 
$V'_r$ each have at least $Q((2^r-1)-(2^{r-1}-2))=Q(2^{r-1}+1) \geq 
3Q$ pebbles.
By Lemma \ref{lem:movable} with $E=2Q$ we can remove $2Q$ pebbles from $V'_1$ and 
place $Q$ of them on $V'_n$. Now $V'_1$ has at least $Q(2^{r-1}-1)$ 
pebbles which is enough to $Q$-cover 
$V'_1,\ldots,V'_{r-1}$ by Theorem \ref{thm:G times P} and $V'_r$ has at least $Q(2^{r-1}+1)$ 
pebbles, more than enough to 
cover $V'_r,\ldots,V'_{n-1}$. Thus $\tD$ is $Q$-coverable, which provides a contradiction. 
Therefore 
one of $V'_{r},\ldots,V'_n$ and $V'_{r+1},\ldots,V'_n, V'_1$ must not be 
saturated, thus satisfying the conclusion of the lemma.
\end{proof}

\begin{lemma}\label{lem:conditions}
There exists a labeling of the vertices of $\tCn$
 such that: 
 \begin{enumerate}
     \item $V_1,..,V_r$ is primary and saturated,
     \item there exists $i\leq r+1$ such that $V_i,..,V_{i+r-1}$ is unsaturated, and
     \item there are no support vertices between $V_1$ and $V_i$.
 \end{enumerate}    
    
\end {lemma}

\begin{proof}
By Lemma \ref{lem:one of each general} there is a labeling such that 
the path
$V_1,\ldots,V_r$ is 
saturated, so it must contain a support vertex.  Let 
$V_k$ be a support vertex with minimum index $k$. Then the 
primary set $V_k,\ldots,V_{k+r-1}$ contains at least as many pebbles as 
$V_1,\ldots,V_r$, therefore it is also saturated.  Let $V'_1=V_k$. By 
Lemma \ref{lem:one of each general} there is an $i$ satisfying the 
second condition.

To show that the third property holds we consider the sets 
\begin{align*}
S=&\{V'_i\,|\,V'_i,\ldots,V'_{i+r-1}\text{ is a saturated primary set}\} \text{ and} \\ 
U=&\{V'_j\,|\,V'_j,\ldots,V'_{j+r-1}\text{ is an unsaturated set}\}. 
\end{align*}
Let $V'_{i^*} \in S$ and $V'_{j^*} \in U$ be such that 
$j^*-i^*=\min\{j-i\, | \, V'_j \in U,\, V'_i \in S,\,\,  j>i\}$. 
By the construction above at least one such pair $i, j$ exists and 
satisfies
$j-i\leq r$, therefore $j^*-i^*\leq r$. Consider the primary saturated set 
$V_{i^*},\ldots,V_{i^*+r-1}$ and the unsaturated set $V_{j^*},\ldots,V_{j^*+r-1}$.

Suppose 
$V'_s$ is a source vertex with $j^*<s<i^*$. If the primary set $V'_s,\ldots,V'_{s+r-1}$ is 
saturated, then we should have replaced $i^*$ with $s$ to obtain a 
lower minimum above. If $V'_s,\ldots,V'_{s+r-1}$ is unsaturated we 
should have
replaced $j^*$ with $s$, again giving a lower minimum. Thus no support 
vertices can lie between $V'_{i^*}$ and $V'_{j^*}$. Finally, relabel the 
vertices so that $V'_{i^*}=V''_1$ to obtain the labeling guaranteed by the 
lemma. (In fact, $i=2$, but we will not be 
using this fact.)
\end{proof}  

Recall that $P=2^r+2^{n-r+1}-3$, $r=\lceil n/2 \rceil$ and we intend to show that 
$\gamma(C_n\square G)\leq P \gamma(G)$.
\begin{theorem} \label{thm:G times C}
    
    Given positive integers $g <Q$, if $\tD$ is any 
       $g$-distribution of $QP$ pebbles on the vertices of $\tCn$, 
       then $\tD$ is $Q$-coverable.

\end{theorem}

\begin{proof}
 
 In search of contradiction suppose $\tD$ is a $g$-distribution on $\tCn$ 
 that is not $Q$-coverable.
  By Lemma  \ref{lem:conditions} we can label the vertices of $\tCn$ so that 
$V_1,...,V_r$ is primary and saturated, $V_i,\ldots,V_{i+r-1}$ is unsaturated, $i \leq r+1$, 
and there are no pebbles on any vertex $V_s$ 
for $1 <s<i$. 
  As there are no support vertices between $V_1$ and $V_i$, the pebbles in  
  $V_2,\ldots,V_r$ are also pebbles in $V_i,\ldots,V_{i+r-1}$.  
  However, this was an unsaturated set, so $|V_2,\ldots,V_r|\leq 
  |V_i,\ldots,V_{i+r-1}|\leq 2^{r}-2$. 
  Thus $V_1$ must have at least $(2^r-1)-|V_2,\ldots,V_r|$ pebbles because 
  $V_1,\ldots,V_r$ was chosen to be saturated. Let 
  $a=(2^r-1)-|V_2,\ldots,V_r|$. Then we can write 
  $|V_1|$ as the sum of two integers, $a$ and $b$, so that $|V_2,\ldots,V_r|+a=Q(2^{r}-1)$ 
  and thus $|V_{r+1},\ldots,V_n|+b=Q(2^{n-r+1}-2)$. Use all of the 
  pebbles on $V_2,..,V_r$ and $a$ pebbles from $V_1$ to $Q$-cover
  the path $V_1,\ldots,V_r$. This is possible by Theorem \ref{thm:G 
  times P}. Now there are at least $Q$ pebbles on each of 
  $V_1,\ldots,V_r$ and at least $Q+b$ pebbles on $V_1$. Consider 
  the path $V_{r+1},\ldots,V_n, V_1$ which contains $n-r+1$ vertices. 
 On this path there are at least $|V_{r+1},\ldots,V_n|+b+Q= Q(2^{n-r+1}-1)$ pebbles and 
  it is therefore $Q$-coverable by Theorem \ref{thm:G times P}. Thus $\tD$ is 
  $Q$-coverable, contradicting the assumption. 
  
\end{proof}
\begin{corollary} \label{cor:cycle inequality}
 $\gamma(G\square 
  C_n)\leq\gamma(G)(2^r+2^{n-r+1}-3)$ for any graph $G$.  
\end{corollary}

 \begin{proof}
     Let $D$ be any distribution of $\gamma(G)P$ pebbles on $(C_n \square 
     G)$. Let $\gamma(G)=Q$ and $g$ be the number of vertices in $G$. 
     By Theorem \ref{thm:G times C} we conclude that the associated 
     $g$-distribution $\tD$ on $\tCn$ is $Q$-coverable. By Lemma 
     \ref{lem:equivalent} the distribution $D$ on $(C_n \square 
     G)$ must also be coverable.
     
  \end{proof}   
 \begin{corollary} \label{cor:cycle}
  $\gamma(C_n)= 2^r+2^{n-r+1}-3$ and $C_n$ is good. 
  \end{corollary}
  
  \begin{proof}
      By Corollary \ref{cor:cycle inequality} with $G=P_1$ we need only show  
      $\gamma(C_n) \geq P$. We number the vertices of $C_n$ sequentially. Consider a 
       distribution with all pebbles placed on $v_1$. The distance 
       from $v_1$ to $v_i$ is $i-1$ when $i\leq r$ and 
       $n-i+1$ when $i>r$. So we require 
       $\sum_{i=1}^r 2^{i-1}+\sum_{i=r+1}^{n} 
       2^{n-i+1}=(2^{r}-1)+(2^{n-r+1}-2)=P$ pebbles.

    \end{proof}  
 
\section{Pebbling numbers for certain products}

Recall that a graph $G$ is {\em good} if there is a 
distribution with only one support vertex requiring $\gamma(G)$ 
pebbles to cover pebble $G$. It was previously known that paths, 
trees and complete graphs are good \cite{H2}.
Section 3 establishes that cycles are good.
In Theorem \ref{thm:relationship} we will prove that there is a relationship 
between the cover pebbling version of Graham's conjecture (Conjecture 
\ref{conj:equality}) and good graphs. First 
note the following:

  \begin{prop} \label{prop:there is one}
  If $G$ and $H$ are good then there is a simple distribution on $G \square 
  H$ that requires $\gamma(G) \gamma(H)$ pebbles. In 
 particular, $\gamma (G\square H) \geq  \gamma(G) \gamma(H)$. 
  \end{prop}   
  \begin{proof}
      Let $w_1,\ldots,w_g$ and $v_1,\ldots,v_h$ be the vertices of $G$ 
      and $H$ respectively, and say $w_1$ and $v_1$ are key vertices for $G$ and $H$.
      Then $\gamma(G)=\sum_{i=1}^g 
      2^{\dist(w_1,w_i)}$ and $\gamma(H)=\sum_{j=1}^h 2^{\dist(v_1,v_j)}$. 
      Consider the distribution on $G \square H$ consisting of a single support 
      vertex $(w_1,v_1)$. This distribution requires 
      
      \begin{align*}
      &\sum_{i=1,j=1}^{i=g, j=h} 2^{\dist((w_1, v_1),(w_i,v_j))}
      =\sum_{i=1,j=1}^{i=g, j=h} 2^{\dist(w_1,w_i)+\dist(v_1,v_j)}\\
=&\sum_{i=1}^g 2^{\dist(w_1,w_i)}\sum_{j=1}^h 2^{\dist(v_1,v_j)}
=\gamma(G)\gamma(H)
      \end{align*}
      pebbles.
   \end{proof}   
  
  \begin{theorem} \label{thm:relationship}
   Suppose $G$ and $H$ are good. Then $\gamma (G \square H)=\gamma(G) 
   \gamma(H)$ if and only if $G \square H$ is good.  
  \end{theorem}    
  
  \begin{proof}
     If $\gamma (G \square H)=\gamma(G)\gamma(H)$  then $G 
     \square H$ is good by Proposition \ref{prop:there is one}.
     If $G \square H$ is good, then by the same proposition it follows 
     that $\gamma (G \square H) \geq\gamma(G)  \gamma(H)$.
      
  Any simple distribution on $G \square H$
     supported on $(w, v)$ would require  $$\sum_{i=1,j=1}^{h,g} 2^{\dist((w, v),(w_i,v_j))}=
      \sum_{i=1}^g 
      2^{\dist(w,w_i)}\sum_{j=1}^h 2^{\dist(v,v_j)}\leq \gamma(G) \times 
      \gamma(H)$$ pebbles, thus proving the other direction of the 
      inequality and concluding the proof of the theorem.
   \end{proof}   
      
\begin{lemma}\label{lem:equality}
  If $G$ is a good graph then $\gamma(P_n \square G)=\gamma(P_n) 
   \gamma (G)$ and $\gamma(C_n \square G)=\gamma(C_n)\gamma (G)$. 
  
 \end{lemma}   
 \begin{proof}
  $P_n$ and $C_n$ are good graphs by Corollaries \ref{cor:path} and \ref 
  {cor:cycle}
  respectively. Let $H_n$ indicate $P_n$ or $C_n$. By Proposition \ref{prop:there is one} 
 we know $\gamma (G\square H_n) \geq  \gamma(G) \gamma(H_n)$. By Corollaries
 \ref{cor:path inequality} and \ref{cor:cycle inequality} we have
 $\gamma (G\square H_n) \leq  \gamma(G) \gamma(H_n)$.
 \end{proof}    

 \begin{corollary} \label{cor:products of good}
   The product of any good graph with $P_n$ or $C_n$ is good.  
\end{corollary}  

\begin{proof}
 This is a direct result of Lemma \ref{lem:equality} and Theorem 
 \ref{thm:relationship}.
\end{proof}    

Now we can easily prove the cover pebbling numbers of some families of 
graphs as advertised in the introduction.

For quick reference, we collect all known cover pebbling numbers here.
\begin{itemize}
\item For any tree $T$, $\gamma(T)=max_{v \in V(T)}(\sum_{u \in 
V(T)}2^{\dist(v,u)})$, as in\cite {H2}.
\item $\gamma(K_n)=2n-1$, as in\cite {H2}.
\item $\gamma(P_n)=2^n-1$, as in\cite {H2}.
\item $\gamma(C_n)= 2^r+2^{n-r+1}-3$, where $r=\lceil n/2 \rceil$, by Corollary 
\ref{cor:cycle}.
\end{itemize}

As all of the graphs referenced in the above list are good, we also know 
that the following products are good, with cover pebbling numbers as 
shown below.

\begin{theorem} \label{thm:results} Let $H=(\square_i P_{n_i}) \square 
(\square_j
C_{m_j}) $.
\begin{itemize}
\item $\gamma (H)=\prod_i \gamma (P_{n_i}) \prod_j \gamma(C_{m_j})$\\
In particular,
\begin {itemize}
\item $\gamma(\square_i P_{n_i}) =\prod_i (2^{n_i}-1)$
\item $\gamma(\square_i C_{m_i}) =\prod_i 
(2^{r_{m_i}}+2^{m_i-r_{m_i}+1}-3)$,
where $r_{m_i}=\lceil m_i/2 \rceil$
\end{itemize}
\item $\gamma(H \square T)=\gamma(H)\gamma(T)$ for any tree $T$
\item $\gamma(H \square K_n)=\gamma(H)\gamma(K_n)$
\end{itemize}    

\end{theorem}
\begin{proof}
    In each statement the fact that the product is good follows from 
    Corollary \ref{cor:products of good}. The pebbling number then 
    follows from Theorem \ref{thm:relationship}.
 \end{proof}

 We also recover a result announced by Hurlbert:
 \begin{corollary}
  The cover pebbling number of the k-hypercube is $3^k$, i.e. $\gamma (Q^k)=3^k$.   
  \end{corollary}   

 \begin {proof}
 As $Q^k$ is isomorphic to $\square^{k} P_2$, by Theorem 
 \ref{thm:results} (part 1) we have
 $\gamma({Q^k})=\prod^{k}\gamma(P_2)=\prod^{k}(2^2-1)=3^k$.
 \end{proof}


\end{document}